\documentclass[VANCOUVER,STIX1COL]{WileyNJD-v2}
\newcommand{\VM}[1]{\mathbf{#1}}
\newcommand{\subs}[1]{\mathcal{#1}}
\newcommand{\norm}[1]{\left \| #1 \right \|}
\newcommand{\R}{\mathbb{R}}
\newlength\figW
\usepackage{tikz}
\usepackage{pgfplots}
\usepackage{float}
\usepackage{cool}
\usepackage{subcaption}
\pgfplotsset{compat=1.14}
\usepackage{todonotes}
\usepackage{mathtools}

\makeatletter
\renewcommand*\env@matrix[1][*\c@MaxMatrixCols c]{%
  \hskip -\arraycolsep
  \let\@ifnextchar\new@ifnextchar
  \array{#1}}
\makeatother

\articletype{Article Type}%

\received{ }
\revised{ }
\accepted{ }

\begin{document}

\title{Calculating the minimal/maximal eigenvalue of symmetric parametrized matrices using projection\protect}

\author[1]{Koen Ruymbeek*}

\author[1]{Karl Meerbergen}

\author[1]{Wim Michiels}
%

\address[1]{\orgdiv{Department of Computer Science}, \orgname{KU Leuven}, \country{Belgium}}

\corres{*Koen Ruymbeek, \email{koen.ruymbeek@cs.kuleuven.be}}

\presentaddress{ Celestijnenlaan 200A, 3001 Leuven.}

\abstract[Summary]{In applications of linear algebra including nuclear physics and structural dynamics, there is a need to deal with uncertainty in the matrices. We focus on matrices that depend on a set of parameters $\VM{\omega}$ and we are interested in the minimal eigenvalue of a matrix pencil $\left(\VM{A}, \VM{B}\right)$ with $\VM{A},\VM{B}$ symmetric and $\VM{B}$ positive definite. If $\VM{\omega}$ can be interpreted as the realisation of random variables, one may be interested in  statistical moments of the minimal eigenvalue. In order to obtain statistical moments, we need a fast evaluation of the eigenvalue as a function of $\VM{\omega}$. Since this is costly for large matrices, we are looking for a small parametrized eigenvalue problem whose minimal eigenvalue makes a small error with  the minimal eigenvalue of the large eigenvalue problem. The advantage, in comparison with a global polynomial approximation (on which, e.g., the polynomial chaos approximation relies), is that we do not suffer from the possible non-smoothness of the minimal eigenvalue. The small scale eigenvalue problem is obtained by projection of the large scale problem. Our main contribution is that for constructing the subspace we use multiple eigenvectors as well as  derivatives of eigenvectors. We provide theoretical results and document numerical experiments regarding the beneficial effect of adding multiple eigenvectors and derivatives.}



\maketitle
%

\section{Introduction}

Let $\VM{A}(\VM{\omega}), \VM{B}(\VM{\omega})$ be $n \times n$ symmetric matrices which smoothly depend on a parameter $\VM{\omega}$ in a compact subset $\Omega \subset \R^d$. We impose that $\VM{B}(\VM{\omega})$ is positive definite, for $\VM{\omega} \in \Omega$, so all eigenvalues $\lambda(\VM{\omega})$ of the generalized eigenvalueproblem $\left(\VM{A}(\VM{\omega}), \VM{B}(\VM{\omega})\right)$, 
$$\VM{A}(\VM{\omega}) \VM{x}(\VM{\omega}) = \lambda(\VM{\omega}) \VM{B}(\VM{\omega}) \VM{x}(\VM{\omega}),$$ are real, see, e.g, \cite{Saad2011} Ch. 9. The objective is to calculate in an efficient way, an accurate global approximation of the minimal eigenvalue $\lambda_1(\VM{\omega})$ over the whole parameterspace. Standard eigensolvers like the Lanczos' method \cite{Lanczos1950}, Jacobi-Davidson \cite{G.Sleijpen1996} or LOBPCG \citep{Knyazev2001} can be applied for large-scale problems for some points in the parameterspace but such an approach becomes expensive when considering a large number of sample points. 

This paper focuses on large scale eigenvalue problems. One class of examples is the estimation of the coercivity constant of parametrized elliptic partial differential equations \cite{Rozza2008}. Another is the computation of the inf-sup constant of PDEs from the minimal singular value of $\VM{A}(\VM{\omega})$, which is computed from the minimal eigenvalue of $\VM{A}(\VM{\omega})^T \VM{A}(\VM{\omega})$. The minimal eigenvalue problem also plays a role in the characterization of pseudospectra \cite{Karow2010}, \cite{Sirkovic2018}. Related problems are 
the determination of statistical moments of the minimal eigenvalue when $\omega$ is the realization of random variables and  the minimization or maximization of the $j$th largest eigenvalue of a parameter dependent Hermitian matrix; see \cite{KangalF2018Asmf} and the references therein.


There are several methods for this kind of problem. The first category concerns polynomial-based methods. In cite{Andreev2012} the authors approximate $\lambda_1(\VM{\omega})$ by sparse tensor products of Legendre polynomials and in \cite{Ghanem2007} polynomial chaos is advocated, a technique which is frequently used when the parameters are stochastic. These techniques  have difficulties with the possible lack of smoothness of the minimal eigenvalue. Another approach is the so-called Successive Constraint Method (SCM), see \cite{Sirkovic2016} and the references therein. In paper \cite{Sirkovic2016}, lower and upper bounds are calculated for each sample in a finite subset of the parameterspace. The main disadvantage of this method is that it does not provide a way to calculate the minimal eigenvalue of a parameter sample other than the samples in the initial subset. 


The idea of this paper is to approximate the minimal eigenvalue of a large scale matrix pencil by the minimal eigenvalue of a smaller matrix, which we call the reduced problem. The reduced problem is obtained by projection of the large scale matrix on a well-chosen subspace. The reduced problem adopts the same smoothness properties or lack thereof as the original large scale problem, which allows us to use less samples than methods based on smooth approximations such as polynomials.
The novelty of this paper is in the construction of the subspace to build the reduced eigenvalue problem. The subspace is constructed from samples of the associated eigenvector of the large scale matrix in the parameterspace. The selection of sample points is based on a greedy method \cite{Quarteroni2016}, with the aim to minimize the residual norm of the large scale eigenvalue problem with a minimum amount of samples. For each sampling point, a large scale eigenvalue problem is solved. Eigenvalue solvers such as Lanczos' method or Jacobi-Davidson usually compute more information than just one eigenvector approximation. The goal is to explore whether other (freely available) information can be used with the aim to further reduce the number of sample points and, consequently, the computational cost. In particular, we will study higher order Hermite interpolation sampling by inserting partial derivatives of the eigenvector towards the parameters in the subspace as well, and the addition of more than one eigenvector approximation for each sample point. Calculating the maximal eigenvalue or an eigenvalue that meets another condition (second maximal/minimal eigenvalue) can be done in a similar way as the method proposed here.


\medskip
The plan of the paper is as follows. In section 2 we introduce the necessary notation and we give some motivating examples. We explain theoretically why we add the first eigenvectors as well as the partial derivatives of the first eigenvector in section 3. Section 4 gives insight how eigenvectors change over the parameters. Subsequently we explain how we build up our subspace and we present the algorithm. In section 6 we illustrate the theory and algorithms numerically. We finish this paper with some concluding remarks and research possibilities.

\section{Notation and motivating examples} 
 We denote vectors by small bold letters and matrices by large bold letters. A subspace of $\R^n$ is denoted by calligraphic letters. We already mentioned that we restrict ourselves to the case where both $\VM{A}$ and $\VM{B}$ are symmetric, $\VM{B}$ is positive definite; $\VM{A}$ and $\VM{B}$ depend analytically on the parameter for all $\VM{\omega} \in \Omega$. From this last requirement, it follows that the inner product defined by 
$$\left( \VM{x}, \VM{y}\right)_\VM{B} := \VM{x}^T \VM{B} \VM{y}, \VM{x},\VM{y} \in \R^n$$ and the associated $\VM{B}$-norm 
$$\norm{\VM{x}}_\VM{B} := \left( \VM{x}, \VM{x}\right)_\VM{B}, \VM{x} \in \R^n$$ are well-defined. We call $\left(\lambda_i, \VM{x}_i\right), i = 1, \hdots, n$ an \emph{eigenpair} of the generalised eigenvalue problem $\left( \VM{A}, \VM{B} \right)$ if $$\VM{A} \VM{x}_i = \lambda_i \VM{B} \VM{x}_i, \quad \VM{x}_i \neq \VM{0}$$
where $\lambda_i$ and $\VM{x}_i$ are respectively called an \emph{eigenvalue} and an associated \emph{eigenvector}. In this case it is proven in \cite{Parlett1998} that the eigenvalues are real and we assume further that 
\begin{equation}\label{eqn:ordering}
\lambda_1 \leq \lambda_2 \leq \hdots \leq \lambda_n. 
\end{equation} 
Eigenvectors are always assumed to be $\VM{B}$-orthonormal, meaning that $\VM{x}_i \VM{B} \VM{x}_j = \delta_{i,j}, i,j = 1,\hdots, n$. 

We further define \begin{equation}
\VM{\VM{\Lambda}} = \mbox{diag}\left( \lambda_1, \lambda_2, \hdots, \lambda_n \right), \quad \VM{X}= [\VM{x}_1, \VM{x}_2, \hdots, \VM{x}_n] \label{eqn:0}
\end{equation} with $\VM{X}^T \VM{B} \VM{X} = \VM{I}_n$ then we can prove that 
\begin{equation}
\VM{X}^T \VM{A} \VM{X} = \VM{\VM{\Lambda}}.  \label{eqn:SVD_A}
\end{equation}
It immediately follows that if $\VM{A}$ is also positive definite, the eigenvalues of the couple $\left( \VM{A}, \VM{B} \right)$ are strictly positive.
Furthermore, if the eigenvalue is simple for $\VM{\omega} \in \Omega \subset \R^d$, there is an environment around $\VM{\omega}$ where this eigenvalue is differentiable as we enumerate the eigenvalues for all parametervalues in increasing order. If the eigenvalue has a multiplicity $m > 1$ for a given $\VM{\omega}$ and is simple in an open set around $\VM{\omega}$, then it can be decomposed in this open set in each direction into $m$ differentiable curves and something similar can be done for $\VM{B}$-normalized eigenvectors, see e.g \cite{Lax2007} ch. 9. This does not mean that it can also be decomposed in $d$-dimensional surfaces, see Example \ref{vbd:1}.   

 Remark that in view of \eqref{eqn:ordering}, $\VM{X}$ is not continuous in the parameter since in an open set around a point where not all eigenvalues are simple, the order of the eigenvectors may change.


\begin{eexample} \label{vbd:1}
Let $\Omega = [-0.5, 0.5] \times [-0.5, 0.5]$ be the parameterspace and let the matrices be 
$$\VM{A}(\VM{\omega}) = \VM{W}(\VM{\omega}) \begin{bmatrix}
\omega_1+1 & \omega_2 &  &  & \\ 
\omega_2 & -\omega_1+1 &  &  & \\ 
 &  & 3 &  & \\ 
 &  &  & \ddots & \\ 
 &  &  &  & n
\end{bmatrix} \VM{W}(\VM{\omega})^T, \VM{B}(\VM{\omega}) = \VM{I}_n$$ 
with $\VM{W}(\VM{\omega}) = [\VM{w}_1, \VM{w}_2, \hdots, \VM{w}_n]$ an orthonormal ($n \times n$)-matrix for all $\VM{\omega} \in \Omega$. The eigenvalues are $\lambda_1(\VM{\omega}) = -\sqrt{\omega_1^2 + \omega_2^2}+1$, $\lambda_2(\VM{\omega}) = \sqrt{\omega_1^2 + \omega_2^2} +1$ and the other eigenvalues are equal to $3, \hdots, n$. The two minimal eigenvalues are displayed in Figure \ref{fig:ex1}. We observe that the eigenvalues cannot be decomposed in 2 smooth surfaces.

\begin{figure}[h]
\setlength{\figW}{5cm} 
  \centerline{\input{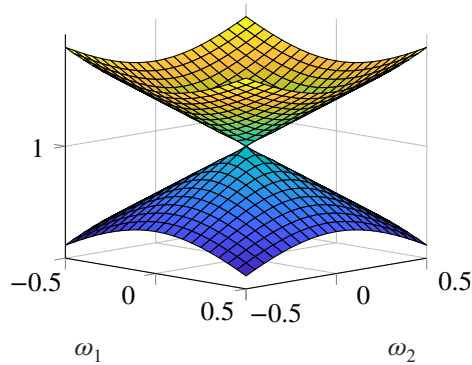}}
\caption{Surfaces that represent $\lambda_1(\VM{\omega})$ and $\lambda_2(\VM{\omega})$ in Example \ref{vbd:1}} \label{fig:ex1}
\end{figure}
\end{eexample}

\begin{eexample}  \label{vbd:2}
This example is taken from \cite[Example 4.4]{Sirkovic2016}. The parameterspace is $\Omega = [0.02, 0.5] \times [2, 8]$ and $\VM{A}(\VM{\omega}) = \sum_{i=1}^3 \theta_i(\VM{\omega}) \VM{A}_i$
with $\theta_i(\VM{\omega}), i = 1,2,3$ analytic functions and $\VM{B}$ is constant. The dimension of the problem is $1311$. The matrices originate from a finite element discretisation of a PDE. In Figure \ref{fig:ex2}, we show the minimal eigenvalue as a function of the parameters. More details can be found in \cite[Example 4.4]{Sirkovic2016}.
\begin{figure}[h]
\begin{center}
\setlength{\figW}{5cm} 
  \centerline{\input{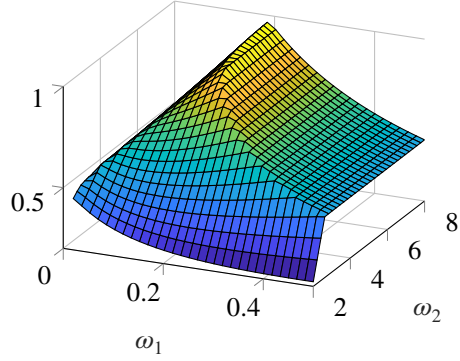}}
\caption{Surface of $\lambda_1(\VM{\omega})$ in Example \ref{vbd:2}} \label{fig:ex2}
\end{center} 
\end{figure} 
\end{eexample}

It is immediately clear from the previous examples that an approximation method (including polynomial approximations) that relies on the smoothness of the surface, is not appropriate. 
We consider projection methods and we describe in the next section why they are well suited in this case.

We end this section by introducing the concept of reduced eigenvalue problems. Let $\subs{V} \subset \R^n$ be a subspace spanned by the columns of an orthonormal matrix $\VM{V} := [\VM{v}_1, \VM{v}_2, \hdots, \VM{v}_m] \in \R^{n\times m}$ with $m$ the dimension of the subspace. We call $\left( \VM{V}^T \VM{A}(\VM{\omega}) \VM{V}, \VM{V}^T \VM{B}(\VM{\omega}) \VM{V} \right)$ the \emph{reduced eigenvalue problem}
on $\subs{V}$. We denote by $\left( \lambda^\subs{V}_i(\VM{\omega}), \VM{x}^\subs{V}_i(\VM{\omega}) \right), i = 1, \hdots, m$ an eigenpair of this reduced eigenvalue problem. 

\section{Hermite interpolation by projection} \label{sect:Herm_int}
In this section, we derive Hermite interpolation properties of the reduced eigenvalue problem when the subspace is built with the first eigenvector and its partial derivatives in some sample points. Concretely, we show that the eigenvalue itself is interpolated as well as its first and second partial derivatives. To prove this, we need a characterisation of the first and second derivative of an eigenvalue and of the first derivative of an eigenvector. Since eigenvectors depend on the normalisation, we first outline the adopted setting, in which they are uniquely defined.
 Let $\lambda_i^*$ be a simple eigenvalue with associated eigenvector $\VM{x}_i^*$ of $\VM{\omega}^* \in \Omega$ then there is an environment $\Omega^*$ around $\VM{\omega}^*$ and smooth functions 
\begin{align}
\begin{aligned}
 \lambda_i:& \Omega^* \rightarrow  \R \\
& \VM{\omega} \mapsto \lambda_i(\VM{\omega})
\end{aligned}
\quad \text{ and } \quad
\begin{aligned}
\VM{x}_i:& \Omega^* \rightarrow  \R^n \\
&\VM{\omega} \mapsto \VM{x}_i(\VM{\omega})
\end{aligned} \label{eqn:setting_deriv}
\end{align} 
such that 
\begin{equation}
\left\{\begin{matrix*}[c]
\VM{A}(\VM{\omega}) \VM{x}_i(\VM{\omega}) - \lambda_i(\VM{\omega}) \VM{B}(\VM{\omega}) \VM{x}_i(\VM{\omega})& = & 0 &,\forall \VM{\omega} \in \Omega^* & \left( (\lambda_i(\VM{\omega}), \VM{x}_i(\VM{\omega})) \text{ is eigenpair} \right)\\ 
\VM{x}_i(\VM{\omega})^T \VM{B}(\VM{\omega}) \VM{x}_i(\VM{\omega})  & = &1 &,\forall \VM{\omega} \in \Omega^* & \left( \text{ normalisation condition} \right) \\ 
\lambda_i(\VM{\omega}^*) & = &\lambda_i^*& & \\
\VM{x}_i(\VM{\omega}^*)  & = &\VM{x}_i^* & &
\end{matrix*}\right. \label{eqn:condit_deriv}
\end{equation}


In the remainder of this section and the next section, we characterise the derivatives of the functions in \eqref{eqn:setting_deriv}. The derivative of a simple eigenvalue $\lambda_i$ with associated eigenvector $\VM{x}_i$ is
\begin{equation}
\pderiv{\lambda_i}{\omega_j} = \VM{x}_i^T \left( \pderiv{\VM{A}}{\omega_j} - \lambda_i \pderiv{\VM{B}}{\omega_j} \right) \VM{x}_i,\quad j = 1, \hdots, d. \label{eqn:dlambda}
\end{equation}
This is a generalisation of the more known result for standard eigenvalue problems, see e.g \cite[Ch. 9]{Lax2007} . For the second derivative, it is sufficient to differentiate equation \eqref{eqn:dlambda}. In this way, we get the following expression for the second derivative
\begin{equation} 
\pderiv[1]{\lambda_i}{\omega_j,\omega_k} = 2 \VM{x}_i^T \left( \pderiv{ \VM{A}}{\omega_j} - \lambda_i \pderiv{\VM{B}}{\omega_j} \right) \pderiv{\VM{x}_i}{\omega_k} + \VM{x}_i^T \left( \pderiv[1]{\VM{A}}{\omega_j, \omega_k} - \pderiv{\lambda_i}{\omega_k} \pderiv{\VM{B}}{\omega_j} - \lambda_i \pderiv[1]{\VM{B}}{\omega_j, \omega_k}  \right) \VM{x}_i, \quad  j,k = 1, \hdots, d. \label{eqn:ddlambda} 
 \nonumber 
\end{equation}
We characterize the derivative of the eigenvector associated with a simple eigenvalue as the solution of a system of linear equations. In this proof we need following lemma.
\begin{lemma} \label{lem:afgel_lambda}
It holds that 
\begin{equation} \label{eqn:help_system_deriv}
\VM{x}_i^T \VM{B} \pderiv{\VM{x}_i}{\omega_j} = -\dfrac{1}{2} \VM{x}_i^T \pderiv{\VM{B}}{\omega_j} \VM{x}_i, \quad j = 1, \hdots, d.
\end{equation}
\begin{proof}
This follows from taking the derivative of the normalisation condition in \eqref{eqn:condit_deriv}.
\end{proof}
\end{lemma}

\begin{theorem} \label{ste:system_deriv}
The derivative of the eigenvector $\VM{x}_i$ associated with the simple eigenvalue $\lambda_i$ is characterized as the solution of the following system:
\begin{equation} \begin{bmatrix}
\lambda_i \VM{B} - \VM{A} & \VM{B} \VM{x}_i \\ 
\VM{x}_i^T \VM{B} & 0
\end{bmatrix} \begin{bmatrix} \pderiv{\VM{x}_i}{\omega_j} \\ \pderiv{\lambda_i}{\omega_j} \end{bmatrix} = \begin{bmatrix}
\left( \pderiv{\VM{A}}{\omega_j} - \lambda_i \pderiv{\VM{B}}{\omega_j} \right ) \VM{x}_i \\
- \dfrac{ \VM{x}_i^T \pderiv{\VM{B}}{\omega_j} \VM{x}_i}{2} 
\end{bmatrix}, \quad j = 1,\hdots, d \label{eqn:system_dx}. 
\end{equation}
\begin{proof}
We start by differentiating $\left( \lambda_i \VM{B} - \VM{A} \right)\VM{x}_i = 0 $ to get 

$$\left( \pderiv{\lambda_i}{\omega_j} \VM{B} + \lambda_i \pderiv{\VM{B}}{\omega_j} - \pderiv{\VM{A}}{\omega_j} \right) \VM{x}_i + \left( \lambda_i \VM{B} - \VM{A} \right) \pderiv{\VM{x}_i}{\omega_i} = 0, \quad j = 1, \hdots, d.$$
As $\left( \lambda_i \VM{B} - \VM{A} \right)$ is a singular matrix, the equation 
\begin{equation} \label{eqn:sys_part1} \left( \lambda_i \VM{B} - \VM{A} \right) \pderiv{\VM{x}_i}{\omega_j} = -\left( \pderiv{\lambda_i}{\omega_j} \VM{B} + \lambda_i \pderiv{\VM{B}}{\omega_j} - \pderiv{\VM{A}}{\omega_j} \right) \VM{x}_i, \quad j = 1, \hdots, d
\end{equation}
is not sufficient to characterise the partial derivative of an eigenvector. The missing information is the information in the direction of $\VM{x}_i$ which is given in \eqref{eqn:help_system_deriv}. The combination of \eqref{eqn:sys_part1} and \eqref{eqn:help_system_deriv} leads to the proof of the theorem. 
\end{proof}
\end{theorem}

We prove in the next property that if both the eigenvector and its partial derivatives are present in the subspace $\subs{V}$ for $\VM{\omega}^* \in \Omega$, the minimal eigenvalue of the reduced eigenvalue problem on $\subs{V}$ is a Hermite interpolant of degree two of the minimal eigenvalue of the pencil $\left( \VM{A}(\VM{\omega}^*), \VM{B}(\VM{\omega}^*) \right)$. Furthermore we state also the well-known property that the minimal eigenvalue of a reduced eigenvalue problem is bounded from below by the original eigenvalue problem.
 
\begin{property} \label{prop:interpol}
  
\begin{enumerate}
\item (Hermite interpolation) If $\left( \lambda_i, \VM{x}_i \right)$ is a simple eigenpair and if $\VM{x}_i \in \subs{V}$, then 
\begin{itemize}
\item $\left(\lambda_i^\subs{V}, \VM{x}_i^\subs{V} \right) = \left(\lambda_i, \VM{V}^T \VM{x}_i \right)$ is an eigenpair of $ \left( \VM{V}^T \VM{A}\VM{V}, \VM{V}^T \VM{B} \VM{V} \right)$.
\item $ \pderiv{\lambda_i^\subs{V}}{\omega_j} = \pderiv{\lambda_i}{\omega_j}, \quad j = 1, \hdots d$.
\end{itemize} 
If $\VM{x}_i, \pderiv{\VM{x}_i}{\omega_j} \in \subs{V}$, then 
$$\pderiv{\lambda_i^\subs{V}}{\omega_j, \omega_k} = \pderiv{\lambda_i}{\omega_j, \omega_k}, \quad k = 1,\hdots, d$$ 
\item If $\subs{V}_1 \subset \subs{V}_2$ then it holds that 
$$\lambda_1 \leq \lambda_1^{\subs{V}_2} \leq \lambda_1^{\subs{V}_1}, \quad \forall \omega \in \Omega, \quad j = 1,\hdots, d.$$
This means that if we extend the subspace, the estimation will be at least as good.
\end{enumerate}
\begin{proof}
\underline{Assertion 1:} \\
The first two statements are well-known results in the case of a standard eigenvalue problem, see e.g. \citep{KangalF2018Asmf} and the references therein. The proofs for the generalized eigenvalue problem are stated for the sake of completeness. 
For the first statement, it is sufficient to see that $\VM{x}_i = \VM{V} \VM{V}^T \VM{x}_i$ to obtain $\VM{V}^T \VM{A} \VM{V}  \VM{V}^T \VM{x}_i = \lambda_i \VM{V}^T \VM{B} \VM{V} \VM{V}^T \VM{x}_i$. \\
For the second, the result follows from 
\begin{align*}
\pderiv{\lambda_i}{\omega_j} & = \VM{x}_i^T \left( \pderiv{\VM{A}}{\omega_j} - \lambda_i \pderiv{\VM{B}}{\omega_j} \right) \VM{x}_i \\
& = \left( \VM{V} \VM{x}_i^\subs{V} \right)^T \left( \pderiv{\VM{A}}{\omega_j} - \lambda_i \pderiv{\VM{B}}{\omega_j} \right) \VM{V} \VM{x}_i^\subs{V} \\
& = \left( \VM{x}_i^\subs{V}\right)^T  \left( \VM{V}^T \pderiv{\VM{A}}{\omega_j}\VM{V} - \lambda_i^\subs{V} \VM{V}^T \pderiv{\VM{B}}{\omega_j} \VM{V} \right) \VM{x}_i^\subs{V} \\ 
& = \pderiv{\lambda_i^\subs{V}}{\omega_j}.
\end{align*}
If also $\pderiv{\VM{x}_i}{\omega_j} \in \subs{V}$, then there exists a vector $\VM{z}$ such that $\pderiv{\VM{x}_i}{\omega_j} = \VM{V} \VM{z}$. We prove that $\VM{z} = \pderiv{\VM{x}_i^\subs{V}}{\omega_j}$, i.e. the projection of $\pderiv{\VM{x}_i}{\omega_j}$ on $\subs{V}$ equals the derivative of the first eigenvector of the reduced eigenvalueproblem.
This derivative is uniquely determined by the system in \eqref{eqn:system_dx}, from which it follows that
$$
\left \{\begin{matrix}
\left(\lambda_i\VM{B} - \VM{A}\right) \VM{V} \VM{z} + \pderiv{\lambda^\subs{V}_i}{\omega_j} \VM{B} \VM{V} \VM{x}_i^\subs{V}& = & \left( \pderiv{\VM{A}}{\omega_j} - \lambda^\subs{V}_i \pderiv{\VM{B}}{\omega_j} \right) \VM{V} \VM{x}_i^\subs{V} \\ 
\left(\VM{x}_i^\subs{V}\right)^T \VM{V}^T \VM{B} \VM{z} & = & - \dfrac{ \left( \VM{x}_i, \VM{V}^T \pderiv{\VM{B}}{\omega_j}  \VM{V} \VM{x}_i\right)}{2} 
\end{matrix} \right.
$$
which means that $[\VM{z}, \pderiv{\lambda^\subs{V}_i}{\omega_j}]^T$ is the vector such that 
$$
\left \{\begin{matrix}
\VM{V}^T \left(\lambda_i \VM{B} - \VM{A} \right) \VM{V} \VM{z} + \pderiv{\lambda^\subs{V}_i}{\omega_j} \VM{V}^T \VM{B} \VM{V} \VM{x}_i^\subs{V}& = & \VM{V}^T \left( \pderiv{\VM{A}}{\omega_j} - \lambda^\subs{V}_i \pderiv{\VM{B}}{\omega_j} \right) \VM{V} \VM{x}_i^\subs{V} \\ 
\left(\VM{x}_i^\subs{V}\right)^T \VM{V}^T \VM{B} \VM{V} \VM{z} & = & - \dfrac{ \left( \VM{x}_i, \VM{V}^T \pderiv{\VM{B}}{\omega_j}  \VM{V} \VM{x}_i\right)}{2} 
\end{matrix} \right. .
$$
This is the system that uniquely determines the derivative of the first eigenvector of the reduced eigenvalue problem, so $\VM{z} = \pderiv{\VM{x}_i^\subs{V}}{\omega_j}$.
From equation \eqref{eqn:ddlambda}, we obtain 
\begin{align*}
\pderiv{\lambda_i}{\omega_j, \omega_k} & = 2 \left( \VM{x}_i^\subs{V}\right)^T \VM{V}^T \left( \pderiv{ \VM{A}}{\omega_j} - \lambda_i \pderiv{\VM{B}}{\omega_j} \right) \VM{V} \pderiv{\VM{x}_i}{\omega_k}^\subs{V} +  \hdots \\
& \quad \left( \VM{x}_i^\subs{V}\right)^T \VM{V}^T \left( \pderiv[1]{\VM{A}}{\omega_j, \omega_k} - \pderiv{\lambda_i}{\omega_k} \pderiv{\VM{B}}{\omega_j} - \lambda_i \pderiv[1]{\VM{B}}{\omega_j, \omega_k}  \right) \VM{V} \VM{x}_i^\subs{V}. \\
& = 2 \left( \VM{x}_i^\subs{V}\right)^T \left( \VM{V}^T \pderiv{ \VM{A}}{\omega_j} \VM{V} - \lambda_i^\subs{V} \VM{V}^T \pderiv{\VM{B}}{\omega_j} \VM{V} \right)  \pderiv{\VM{x}_i^\subs{V}}{\omega_k} +  \hdots \\
& \quad \left( \VM{x}_i^\subs{V}\right)^T  \left( \VM{V}^T \pderiv[1]{\VM{A}}{\omega_j, \omega_k} \VM{V} - \pderiv{\lambda_i^\subs{V}}{\omega_k} \VM{V}^T \pderiv{\VM{B}}{\omega_j}\VM{V} - \lambda_i^\subs{V} \VM{V}^T \pderiv[1]{\VM{B}}{\omega_j, \omega_k} \VM{V} \right) \VM{x}_i^\subs{V} \\
& = \pderiv{\lambda_i^\subs{V}}{\omega_j, \omega_k}.
\end{align*}
\underline{Assertion 2:} This results is well known, we refer to \cite{Parlett1998} for a proof.
\end{proof}
\end{property}

\section{Characterisation of the derivative of an eigenvector}
The previous section showed that adding the derivative of the eigenvector to the subspace leads to higher-order Hermite interpolation in the eigenvalue. The aim of this section is to get more insight in the behaviour of eigenvalues as a function of the parameters by deriving an analytic formula for the derivative of the eigenvector. The next result can also be found in \cite{seyranian2003multiparameter}, but we give here our proof, which is based on the diagonalisation of $\lambda_1 B-A$ in \eqref{eqn:system_dx}, to show the connection with Theorem \ref{ste:system_deriv}.


\begin{theorem} \label{ste:dXdw}
If $\VM{x}_i$ is an eigenvector associated with a simple eigenvalue $\lambda_i$, then we have
\begin{equation} \label{eqn:dxdomega}
\pderiv{\VM{x}_i}{\omega_j} = - \dfrac{\VM{x}_i^T  \pderiv{\VM{B}}{\omega_j} \VM{x}_i }{2} \VM{x}_i  + \sum_{k=1, k \neq i}^n \left( \VM{x}_k^T \dfrac{ \left( \pderiv{\VM{A}}{\omega_j} - \lambda_i \pderiv{\VM{B}}{\omega_j} \right) \VM{x}_i}{\lambda_i - \lambda_k} \right) \VM{x}_k, \quad j = 1, \hdots, d.
\end{equation} 

\begin{proof}
Without loss of generality, we prove the statement for $i = 1$.
From \eqref{eqn:SVD_A} it follows immediately that system \eqref{eqn:system_dx} can be written as
$$\begin{bmatrix}[c|c]
\lambda_1 \left(\VM{X}^{-1} \right)^T \VM{X}^{-1} - \left(\VM{X}^{-1} \right)^T \VM{\Lambda} \VM{X}^{-1} & \left(\VM{X}^{-1} \right)^T \VM{e}_1 \\ \hline 
\VM{e}_1^T  \VM{X}^{-1}  & 0
\end{bmatrix} \begin{bmatrix} \pderiv{\VM{x}_1}{\omega_j} \\ \pderiv{\lambda_1}{\omega_j} \end{bmatrix} = \begin{bmatrix}
\left( \pderiv{\VM{A}}{\omega_j} - \lambda_1 \pderiv{\VM{B}}{\omega_j} \right ) \VM{x}_1 \\ \hline
- \dfrac{  \VM{x}_1^T \pderiv{\VM{B}}{\omega_j} \VM{x}_1}{2} 
\end{bmatrix}. $$
Using $\VM{X} \VM{X}^T \VM{B} = \VM{I} $, we obtain
$$\begin{bmatrix}[c|c]
\lambda_1 \left(\VM{X}^{-1} \right)^T  - \left(\VM{X}^{-1} \right)^T \VM{\Lambda} & \left(\VM{X}^{-1} \right)^T \VM{e}_1 \\ \hline
\VM{e}_1^T    & 0
\end{bmatrix} \begin{bmatrix} \VM{X} ^T \VM{B} \pderiv{\VM{x}_1}{\omega_j} \\ \pderiv{\lambda_1}{\omega_j} \end{bmatrix} = \begin{bmatrix}
\left( \pderiv{\VM{A}}{\omega_j} - \lambda_1 \pderiv{\VM{B}}{\omega_j} \right ) \VM{x}_1 \\ \hline
- \dfrac{  \VM{x}_1^T \pderiv{\VM{B}}{\omega_j} \VM{x}_1}{2} 
\end{bmatrix}. $$
By multiplying from the left with the non-singular matrix  $\begin{bmatrix}
\VM{X}^T & \VM{0} \\ 0 & 1
\end{bmatrix}$ we get
$$\begin{bmatrix}[c|c]
\lambda_1 \VM{I} - \VM{\Lambda} & \VM{e}_1 \\ \hline \VM{e}_1^T & 0 
\end{bmatrix} \begin{bmatrix} \VM{X} ^T \VM{B} \pderiv{\VM{x}_1}{\omega_i} \\ \pderiv{\lambda_1}{\omega_j} \end{bmatrix} = \begin{bmatrix}
\VM{X}^T \left( \pderiv{\VM{A}}{\omega_j} - \lambda_1 \pderiv{\VM{B}}{\omega_j} \right) \VM{x}_1 \\ \hline
- \dfrac{  \VM{x}_1^T \pderiv{\VM{B}}{\omega_j} \VM{x}_1}{2}\end{bmatrix}.
$$
From this expression it follows that

\begin{align*}
\pderiv{\lambda_1}{\omega_j} & = \VM{x}_1^T \left( \pderiv{\VM{A}}{\omega_j} - \lambda_1 \pderiv{\VM{B}}{\omega_j} \right) \VM{x}_1 \\
\VM{x}_k^T \VM{B} \pderiv{\VM{x}_1}{\omega_j} & = \left\{\begin{matrix}
- \dfrac{  \VM{x}_k^T \pderiv{\VM{B}}{\omega_j} \VM{x}_1}{2} & ,k = 1\\ 
\dfrac{\VM{x}_k^T \left( \pderiv{\VM{A}}{\omega_j} - \lambda_1 \pderiv{\VM{B}}{\omega_j} \right) \VM{x}_1}{\lambda_1 - \lambda_k} & ,k \neq 1.
\end{matrix}\right.
\end{align*}
Finally, multiplying with $\VM{X}$ from the left proves the statement.
\end{proof}
\end{theorem}

We observe that the weight of $\VM{x}_k, k = 2, \hdots, n$ in the expression for $\pderiv{\VM{x}_1}{\omega_j}$ depends on how close $\lambda_k$ is to $\lambda_1$ and depends on the projection of  $\left( \pderiv{\VM{A}}{\omega_j} - \lambda_1 \pderiv{\VM{B}}{\omega_j} \right) \VM{x}_1, j = 1, \hdots,d$ on $\VM{x}_k$. In general, the closer $\lambda_k$ is to $\lambda_1$ the higher is the impact of $\VM{x}_k$. We also deduce that if $\VM{B}$ does not depend on the variables then the partial derivative is $\VM{B}$-orthonormal to $\VM{x}_1$.
Example \ref{ex:deriv} shows that the norm of the partial derivative of the eigenvector cannot be bounded uniformly on the set $\Omega$, despite $\VM{A}$ and $\VM{B}$ being analytic and $\VM{B}$ positive definite for $\VM{\omega} \in \Omega$.

\begin{eexample} \label{ex:deriv}
We reconsider example \ref{vbd:1} with $n = 2$ and $\VM{W} = \VM{I}_2$. In Figure \ref{fig:example_3}, we plot the eigenvalues $\lambda_1(\VM{\omega})$ and $\lambda_2(\VM{\omega})$ in (a) and (b) and the norm of $\pderiv{\VM{x}_1}{\omega_2}$ in (c) and (d) as a function of $\omega_2$, for fixed values of $\omega_1$.
\begin{figure}[h]
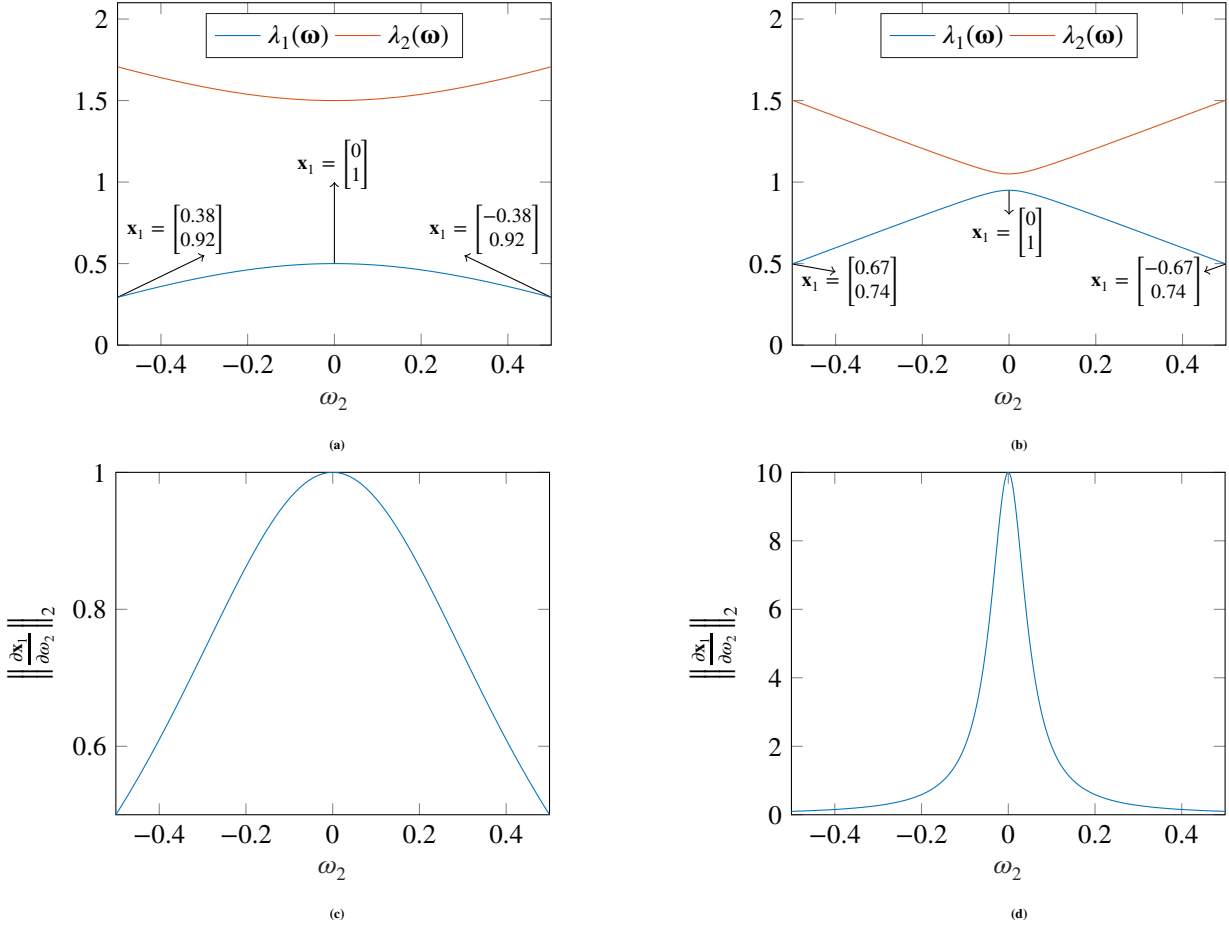

\setlength{\figW}{6cm} 
 \begin{subfigure}[b]{0.5\textwidth}
  \hspace{0.8cm}
%
%
%
\definecolor{mycolor1}{rgb}{0.00000,0.44700,0.74100}%
\definecolor{mycolor2}{rgb}{0.85000,0.32500,0.09800}%
\begin{tikzpicture}

\begin{axis}[%
width=0.951\figW,
height=0.75\figW,
at={(0\figW,0\figW)},
scale only axis,
xmin=-0.500000,
xmax=0.500000,
xlabel style={font=\color{white!15!black}},
xlabel={$\omega_2$},
ymin=0.000000,
ymax=2.100000,
axis background/.style={fill=white},
legend style={at={(0.5,0.97)}, anchor=north, legend columns=2, legend cell align=left, align=left, draw=white!15!black}
]
\draw [->] (-0.5,0.292893) -- (-0.3,0.55);
\draw [->] (0.5,0.292893) -- (0.3,0.55);
\draw [->] (0, 0.5) -- (0, 1);
\node[above right] at (-0.5,0.5) { \footnotesize $\VM{x}_1 =\begin{bmatrix}
0.38\\ 
0.92
\end{bmatrix}$ };
\node[above left] at (0.5,0.5) { \footnotesize $\VM{x}_1 =\begin{bmatrix}
-0.38\\ 
0.92
\end{bmatrix}$ };
\node[above] at (0, 0.9) { \footnotesize $\VM{x}_1 =\begin{bmatrix}
0\\ 
1
\end{bmatrix}$ };

\addplot [color=mycolor1]
  table[row sep=crcr]{%
-0.500000	0.292893\\
-0.490000	0.299929\\
-0.480000	0.306891\\
-0.470000	0.313778\\
-0.460000	0.320588\\
-0.450000	0.327319\\
-0.440000	0.333967\\
-0.430000	0.340531\\
-0.420000	0.347007\\
-0.410000	0.353393\\
-0.400000	0.359688\\
-0.390000	0.365886\\
-0.380000	0.371987\\
-0.370000	0.377987\\
-0.360000	0.383883\\
-0.350000	0.389672\\
-0.340000	0.395351\\
-0.330000	0.400917\\
-0.320000	0.406367\\
-0.310000	0.411697\\
-0.300000	0.416905\\
-0.290000	0.421986\\
-0.280000	0.426938\\
-0.270000	0.431757\\
-0.260000	0.436440\\
-0.250000	0.440983\\
-0.240000	0.445383\\
-0.230000	0.449636\\
-0.220000	0.453740\\
-0.210000	0.457690\\
-0.200000	0.461484\\
-0.190000	0.465117\\
-0.180000	0.468587\\
-0.170000	0.471890\\
-0.160000	0.475024\\
-0.150000	0.477985\\
-0.140000	0.480770\\
-0.130000	0.483376\\
-0.120000	0.485802\\
-0.110000	0.488043\\
-0.100000	0.490098\\
-0.090000	0.491965\\
-0.080000	0.493640\\
-0.070000	0.495124\\
-0.060000	0.496413\\
-0.050000	0.497506\\
-0.040000	0.498403\\
-0.030000	0.499101\\
-0.020000	0.499600\\
-0.010000	0.499900\\
0.000000	0.500000\\
0.010000	0.499900\\
0.020000	0.499600\\
0.030000	0.499101\\
0.040000	0.498403\\
0.050000	0.497506\\
0.060000	0.496413\\
0.070000	0.495124\\
0.080000	0.493640\\
0.090000	0.491965\\
0.100000	0.490098\\
0.110000	0.488043\\
0.120000	0.485802\\
0.130000	0.483376\\
0.140000	0.480770\\
0.150000	0.477985\\
0.160000	0.475024\\
0.170000	0.471890\\
0.180000	0.468587\\
0.190000	0.465117\\
0.200000	0.461484\\
0.210000	0.457690\\
0.220000	0.453740\\
0.230000	0.449636\\
0.240000	0.445383\\
0.250000	0.440983\\
0.260000	0.436440\\
0.270000	0.431757\\
0.280000	0.426938\\
0.290000	0.421986\\
0.300000	0.416905\\
0.310000	0.411697\\
0.320000	0.406367\\
0.330000	0.400917\\
0.340000	0.395351\\
0.350000	0.389672\\
0.360000	0.383883\\
0.370000	0.377987\\
0.380000	0.371987\\
0.390000	0.365886\\
0.400000	0.359688\\
0.410000	0.353393\\
0.420000	0.347007\\
0.430000	0.340531\\
0.440000	0.333967\\
0.450000	0.327319\\
0.460000	0.320588\\
0.470000	0.313778\\
0.480000	0.306891\\
0.490000	0.299929\\
0.500000	0.292893\\
};
\addlegendentry{$\lambda_1(\VM{\omega})$}

\addplot [color=mycolor2]
  table[row sep=crcr]{%
-0.500000	1.707107\\
-0.490000	1.700071\\
-0.480000	1.693109\\
-0.470000	1.686222\\
-0.460000	1.679412\\
-0.450000	1.672681\\
-0.440000	1.666033\\
-0.430000	1.659469\\
-0.420000	1.652993\\
-0.410000	1.646607\\
-0.400000	1.640312\\
-0.390000	1.634114\\
-0.380000	1.628013\\
-0.370000	1.622013\\
-0.360000	1.616117\\
-0.350000	1.610328\\
-0.340000	1.604649\\
-0.330000	1.599083\\
-0.320000	1.593633\\
-0.310000	1.588303\\
-0.300000	1.583095\\
-0.290000	1.578014\\
-0.280000	1.573062\\
-0.270000	1.568243\\
-0.260000	1.563560\\
-0.250000	1.559017\\
-0.240000	1.554617\\
-0.230000	1.550364\\
-0.220000	1.546260\\
-0.210000	1.542310\\
-0.200000	1.538516\\
-0.190000	1.534883\\
-0.180000	1.531413\\
-0.170000	1.528110\\
-0.160000	1.524976\\
-0.150000	1.522015\\
-0.140000	1.519230\\
-0.130000	1.516624\\
-0.120000	1.514198\\
-0.110000	1.511957\\
-0.100000	1.509902\\
-0.090000	1.508035\\
-0.080000	1.506360\\
-0.070000	1.504876\\
-0.060000	1.503587\\
-0.050000	1.502494\\
-0.040000	1.501597\\
-0.030000	1.500899\\
-0.020000	1.500400\\
-0.010000	1.500100\\
0.000000	1.500000\\
0.010000	1.500100\\
0.020000	1.500400\\
0.030000	1.500899\\
0.040000	1.501597\\
0.050000	1.502494\\
0.060000	1.503587\\
0.070000	1.504876\\
0.080000	1.506360\\
0.090000	1.508035\\
0.100000	1.509902\\
0.110000	1.511957\\
0.120000	1.514198\\
0.130000	1.516624\\
0.140000	1.519230\\
0.150000	1.522015\\
0.160000	1.524976\\
0.170000	1.528110\\
0.180000	1.531413\\
0.190000	1.534883\\
0.200000	1.538516\\
0.210000	1.542310\\
0.220000	1.546260\\
0.230000	1.550364\\
0.240000	1.554617\\
0.250000	1.559017\\
0.260000	1.563560\\
0.270000	1.568243\\
0.280000	1.573062\\
0.290000	1.578014\\
0.300000	1.583095\\
0.310000	1.588303\\
0.320000	1.593633\\
0.330000	1.599083\\
0.340000	1.604649\\
0.350000	1.610328\\
0.360000	1.616117\\
0.370000	1.622013\\
0.380000	1.628013\\
0.390000	1.634114\\
0.400000	1.640312\\
0.410000	1.646607\\
0.420000	1.652993\\
0.430000	1.659469\\
0.440000	1.666033\\
0.450000	1.672681\\
0.460000	1.679412\\
0.470000	1.686222\\
0.480000	1.693109\\
0.490000	1.700071\\
0.500000	1.707107\\
};
\addlegendentry{$\lambda_2(\VM{\omega})$}

\end{axis}
\end{tikzpicture}%
  \caption{}
\end{subfigure}
\begin{subfigure}[b]{0.5\textwidth}
	\hspace{0.7cm}  
%
%
%
\definecolor{mycolor1}{rgb}{0.00000,0.44700,0.74100}%
\definecolor{mycolor2}{rgb}{0.85000,0.32500,0.09800}%
\begin{tikzpicture}

\begin{axis}[%
width=0.951\figW,
height=0.75\figW,
at={(0\figW,0\figW)},
scale only axis,
xmin=-0.500000,
xmax=0.500000,
xlabel style={font=\color{white!15!black}},
xlabel={$\omega_2$},
ymin=0.000000,
ymax=2.100000,
axis background/.style={fill=white},
legend style={at={(0.5,0.97)}, anchor=north, legend columns=2, legend cell align=left, align=left, draw=white!15!black}
]
\draw [->] (-0.5,0.497506) -- (-0.4,0.45);
\draw [->] (0.5,0.497506) -- (0.45,0.45);
\draw [->] (0, 0.95) -- (0, 0.8);
\node[below right] at (-0.5,0.6) { \footnotesize $\VM{x}_1 =\begin{bmatrix}
0.67\\ 
0.74
\end{bmatrix}$ };
\node[below left] at (0.47,0.6) { \footnotesize $\VM{x}_1 =\begin{bmatrix}
-0.67\\ 
0.74
\end{bmatrix}$ };
\node[below] at (0,0.9) { \footnotesize $\VM{x}_1 =\begin{bmatrix}
0\\ 
1
\end{bmatrix}$ };

\addplot [color=mycolor1]
  table[row sep=crcr]{%
-0.500000	0.497506\\
-0.490000	0.507456\\
-0.480000	0.517403\\
-0.470000	0.527348\\
-0.460000	0.537291\\
-0.450000	0.547231\\
-0.440000	0.557168\\
-0.430000	0.567103\\
-0.420000	0.577034\\
-0.410000	0.586962\\
-0.400000	0.596887\\
-0.390000	0.606808\\
-0.380000	0.616725\\
-0.370000	0.626637\\
-0.360000	0.636544\\
-0.350000	0.646447\\
-0.340000	0.656343\\
-0.330000	0.666234\\
-0.320000	0.676117\\
-0.310000	0.685994\\
-0.300000	0.695862\\
-0.290000	0.705721\\
-0.280000	0.715571\\
-0.270000	0.725409\\
-0.260000	0.735236\\
-0.250000	0.745049\\
-0.240000	0.754847\\
-0.230000	0.764628\\
-0.220000	0.774390\\
-0.210000	0.784130\\
-0.200000	0.793845\\
-0.190000	0.803531\\
-0.180000	0.813185\\
-0.170000	0.822800\\
-0.160000	0.832369\\
-0.150000	0.841886\\
-0.140000	0.851339\\
-0.130000	0.860716\\
-0.120000	0.870000\\
-0.110000	0.879170\\
-0.100000	0.888197\\
-0.090000	0.897044\\
-0.080000	0.905660\\
-0.070000	0.913977\\
-0.060000	0.921898\\
-0.050000	0.929289\\
-0.040000	0.935969\\
-0.030000	0.941690\\
-0.020000	0.946148\\
-0.010000	0.949010\\
0.000000	0.950000\\
0.010000	0.949010\\
0.020000	0.946148\\
0.030000	0.941690\\
0.040000	0.935969\\
0.050000	0.929289\\
0.060000	0.921898\\
0.070000	0.913977\\
0.080000	0.905660\\
0.090000	0.897044\\
0.100000	0.888197\\
0.110000	0.879170\\
0.120000	0.870000\\
0.130000	0.860716\\
0.140000	0.851339\\
0.150000	0.841886\\
0.160000	0.832369\\
0.170000	0.822800\\
0.180000	0.813185\\
0.190000	0.803531\\
0.200000	0.793845\\
0.210000	0.784130\\
0.220000	0.774390\\
0.230000	0.764628\\
0.240000	0.754847\\
0.250000	0.745049\\
0.260000	0.735236\\
0.270000	0.725409\\
0.280000	0.715571\\
0.290000	0.705721\\
0.300000	0.695862\\
0.310000	0.685994\\
0.320000	0.676117\\
0.330000	0.666234\\
0.340000	0.656343\\
0.350000	0.646447\\
0.360000	0.636544\\
0.370000	0.626637\\
0.380000	0.616725\\
0.390000	0.606808\\
0.400000	0.596887\\
0.410000	0.586962\\
0.420000	0.577034\\
0.430000	0.567103\\
0.440000	0.557168\\
0.450000	0.547231\\
0.460000	0.537291\\
0.470000	0.527348\\
0.480000	0.517403\\
0.490000	0.507456\\
0.500000	0.497506\\
};
\addlegendentry{$\lambda_1(\VM{\omega})$}

\addplot [color=mycolor2]
  table[row sep=crcr]{%
-0.500000	1.502494\\
-0.490000	1.492544\\
-0.480000	1.482597\\
-0.470000	1.472652\\
-0.460000	1.462709\\
-0.450000	1.452769\\
-0.440000	1.442832\\
-0.430000	1.432897\\
-0.420000	1.422966\\
-0.410000	1.413038\\
-0.400000	1.403113\\
-0.390000	1.393192\\
-0.380000	1.383275\\
-0.370000	1.373363\\
-0.360000	1.363456\\
-0.350000	1.353553\\
-0.340000	1.343657\\
-0.330000	1.333766\\
-0.320000	1.323883\\
-0.310000	1.314006\\
-0.300000	1.304138\\
-0.290000	1.294279\\
-0.280000	1.284429\\
-0.270000	1.274591\\
-0.260000	1.264764\\
-0.250000	1.254951\\
-0.240000	1.245153\\
-0.230000	1.235372\\
-0.220000	1.225610\\
-0.210000	1.215870\\
-0.200000	1.206155\\
-0.190000	1.196469\\
-0.180000	1.186815\\
-0.170000	1.177200\\
-0.160000	1.167631\\
-0.150000	1.158114\\
-0.140000	1.148661\\
-0.130000	1.139284\\
-0.120000	1.130000\\
-0.110000	1.120830\\
-0.100000	1.111803\\
-0.090000	1.102956\\
-0.080000	1.094340\\
-0.070000	1.086023\\
-0.060000	1.078102\\
-0.050000	1.070711\\
-0.040000	1.064031\\
-0.030000	1.058310\\
-0.020000	1.053852\\
-0.010000	1.050990\\
0.000000	1.050000\\
0.010000	1.050990\\
0.020000	1.053852\\
0.030000	1.058310\\
0.040000	1.064031\\
0.050000	1.070711\\
0.060000	1.078102\\
0.070000	1.086023\\
0.080000	1.094340\\
0.090000	1.102956\\
0.100000	1.111803\\
0.110000	1.120830\\
0.120000	1.130000\\
0.130000	1.139284\\
0.140000	1.148661\\
0.150000	1.158114\\
0.160000	1.167631\\
0.170000	1.177200\\
0.180000	1.186815\\
0.190000	1.196469\\
0.200000	1.206155\\
0.210000	1.215870\\
0.220000	1.225610\\
0.230000	1.235372\\
0.240000	1.245153\\
0.250000	1.254951\\
0.260000	1.264764\\
0.270000	1.274591\\
0.280000	1.284429\\
0.290000	1.294279\\
0.300000	1.304138\\
0.310000	1.314006\\
0.320000	1.323883\\
0.330000	1.333766\\
0.340000	1.343657\\
0.350000	1.353553\\
0.360000	1.363456\\
0.370000	1.373363\\
0.380000	1.383275\\
0.390000	1.393192\\
0.400000	1.403113\\
0.410000	1.413038\\
0.420000	1.422966\\
0.430000	1.432897\\
0.440000	1.442832\\
0.450000	1.452769\\
0.460000	1.462709\\
0.470000	1.472652\\
0.480000	1.482597\\
0.490000	1.492544\\
0.500000	1.502494\\
};
\addlegendentry{$\lambda_2(\VM{\omega})$}

\end{axis}
\end{tikzpicture}%
  \caption{}
\end{subfigure} \\
\begin{subfigure}[b]{0.5\textwidth}
  \input{images/eigv_norm_kegel_0.5.tex}
  \caption{}
\end{subfigure}
\begin{subfigure}[b]{0.5\textwidth}
  \input{images/eigv_norm_kegel_0.05.tex}
  \caption{}
\end{subfigure}
\caption{(a) and (b): $\lambda_1(\VM{\omega})$ and $\lambda_2(\VM{\omega})$ with $\omega_1 = 0.5$ (left) and $\omega_1 = 0.05$ (right), (c) and (d): norm of $\pderiv{\VM{x}_1(\VM{\omega})}{\omega_2}$  with $\omega_1= 0.5$ (left) and $\omega_1= 0.05$ (right)} \label{fig:example_3}
\end{figure}
One can observe that for $\omega_2 = 0$, the norm of the partial derivative to $\omega_2$ becomes larger if $\omega_1$ goes to 0. It can be shown 
by Theorem \ref{ste:dXdw} that this norm is unbounded. For $\omega_2 = 0$, the eigenvectors are equal to $$\VM{x}_1(\omega_1, 0) = \begin{bmatrix}
0 \\ 1
\end{bmatrix} \mbox{ and } \VM{x}_2(\omega_1, 0) = \begin{bmatrix}
1 \\ 0
\end{bmatrix}$$
from which we can deduce that
\begin{align*}
\dfrac{ \partial \VM{x}_1(\VM{\omega})}{\partial \omega_2} & = \VM{x}_2(\omega_1, 0)\dfrac{ \left( \pderiv{\VM{A}(\omega_1,0) }{\omega_2} \VM{x}_1(\omega_1, 0), \VM{x}_2(\omega_1, 0)\right)}{\lambda_1(\omega_1,0) - \lambda_2(\omega_1,0)} \\
& = \VM{x}_2(\omega_1, 0)\dfrac{ \left( \begin{bmatrix}
0 & 1\\  1 & 0 \end{bmatrix}\begin{bmatrix} 0 \\ 1 \end{bmatrix}, \begin{bmatrix}
1 \\ 0 \end{bmatrix} \right)}{\lambda_1(\omega_1,0) - \lambda_2(\omega_1,0)} \\
& = \VM{x}_2(\omega_1, 0)\dfrac{ 1}{-2 \omega_1}. \\ 
\end{align*} 
We conclude that if $\omega_1$ goes to $0$, the norm of the partial derivative goes to infinity.
\end{eexample} 
Although the analytic formula for the partial derivative gives much insight how eigenvectors depend on the other eigenvectors and eigenvalues, it is less suitable to compute derivatives. The reason is that we need all eigenvectors and eigenvalues and that in practice it is far too costly to calculate these. In practice we calculate the partial derivatives by solving system \eqref{eqn:system_dx}. Note that the cost of solving \eqref{eqn:system_dx} for all parameters is not negligible, since it depends on the dimension and the sparsity pattern of the matrices. Therefore we will propose an alternative construction of the subspace that adds approximations to the derivatives in the next section. 


\section{Computing a global approximation} \label{sect:glob_approx}
\subsection{Subspace expansion}
Property \ref{prop:interpol} presents the key theory for the construction of the reduced problem: adding eigenvectors and partial derivatives leads to Hermite interpolation of degree two for the eigenvalue. The eigenvalue and associated eigenvector are computed by a subspace projection method. The computation of the derivatives of the eigenvector from \eqref{eqn:system_dx} is usually computationally expensive. Therefore we look for cheaper alternatives. Eq. \eqref{eqn:dxdomega} suggests that the derivatives mainly depend on the eigenvectors associated with eigenvalues closest to $\lambda_1$.In particular, this is the case when the spectrum is well-separated with only a few eigenvalues near $\lambda_1$. This implies that adding a few eigenvectors to the subspace in an interpolation point may serve as a computationally cheap surrogate for adding one eigenvector and its derivatives. Approximations to the associated eigenvectors of the eigenvalues nearest $\lambda_1$ are usually available from subspace methods. Indeed, when we use the Lanczos' method, the subspace produces reasonably good approximations to these eigenvectors, at no additional cost. Moreover, it should be noted that there is another reason why adding the second and third eigenvector may be beneficial: the eigenvector of the 2nd minimal eigenvalue may be a good approximation to the eigenvector associated with the minimal eigenvalue at another parameter value $\VM{\omega}_1 \in \Omega$. This means we probably do not need to solve a costly eigenvalue problem in $\VM{\omega}_1$. Example \ref{ex:vb3} gives a good illustration of such a situation. In addition, we see in the next section that adding the second eigenvector helps making a sharp upper bound for the error. 
%

\begin{eexample} \label{ex:vb3}
We reconsider example \ref{vbd:1} where we take $n$ arbitrarily large.
If $\omega_2 = 0$, then one can verify that
$$\VM{x}_1(\omega_1,0) = \left\{\begin{matrix}
\VM{w}_1(\omega_1,0) & ,\omega_1 \leq 0\\ 
\VM{w}_2(\omega_1,0) & ,\omega_1 > 0
\end{matrix}\right.,
\quad
\VM{x}_2(\omega_1,0) = \left\{\begin{matrix}
\VM{w}_2(\omega_1,0) & ,\omega_1 \leq 0\\ 
\VM{w}_1(\omega_1,0) & ,\omega_1 > 0
\end{matrix}\right. .$$ Take $\omega^1 < 0$ and $\omega^2 > 0$ two points in the first dimension of $\Omega$. If both $\VM{w}_1$ and $\VM{w}_2$ do not strongly depend on $\VM{\omega}$, this means that $\VM{x}_1(\omega^1,0) \approx \VM{x}_2(\omega^2,0)$. Therefore we conclude that $\subs{V}_1 :=\text{span}  \{ \VM{x}_1(\omega^i,0), i = 1,2 \} \approx \text{span} \{ \VM{x}_i(\omega^1,0), i = 1,2 \} =: \subs{V}_2$, but $\subs{V}_2$ is less costly to calculate than $\subs{V}_1$. This is illustrated in Figure \ref{fig:ex4}.

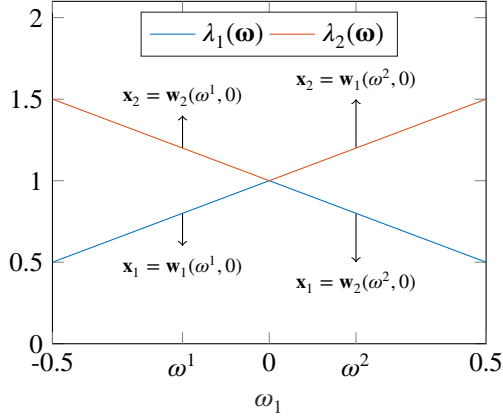
\begin{figure}[h] 
	\begin{center}
%
%
%
\definecolor{mycolor1}{rgb}{0.00000,0.44700,0.74100}%
\definecolor{mycolor2}{rgb}{0.85000,0.32500,0.09800}%
\begin{tikzpicture}

\begin{axis}[%
width=0.951\figW,
height=0.75\figW,
at={(0\figW,0\figW)},
scale only axis,
xmin=-0.500000,
xmax=0.500000,
xtick={\empty},
extra x ticks={-0.5, -0.2, 0, 0.2, 0.5},
extra x tick labels={-0.5, $\omega^1$, 0, $\omega^2$, 0.5},
xlabel style={font=\color{white!15!black}},
xlabel={$\omega_1$},
ymin=0.000000,
ymax=2.100000,
axis background/.style={fill=white},
legend style={at={(0.5,0.97)}, anchor=north, legend columns=2, legend cell align=left, align=left, draw=white!15!black}
]

\draw [->] (-0.2,1.2) -- (-0.2,1.4);
\draw [->] (-0.2,0.8) -- (-0.2,0.6);
\draw [->] (0.2,1.2) -- (0.2,1.5);
\draw [->] (0.2,0.8) -- (0.2,0.5);
\node[above] at (-0.2,1.4) { \footnotesize $\VM{x}_2 = \VM{w}_2(\omega^1,0)$};
\node[above] at (0.2,1.5) { \footnotesize $\VM{x}_2 = \VM{w}_1(\omega^2,0)$};
\node[below] at (-0.2,0.6) { \footnotesize $\VM{x}_1  = \VM{w}_1(\omega^1,0)$};
\node[below] at (0.2,0.5) { \footnotesize $\VM{x}_1 = \VM{w}_2(\omega^2,0) $};

\addplot [color=mycolor1]
  table[row sep=crcr]{%
-0.500000	0.500000\\
-0.490000	0.510000\\
-0.480000	0.520000\\
-0.470000	0.530000\\
-0.460000	0.540000\\
-0.450000	0.550000\\
-0.440000	0.560000\\
-0.430000	0.570000\\
-0.420000	0.580000\\
-0.410000	0.590000\\
-0.400000	0.600000\\
-0.390000	0.610000\\
-0.380000	0.620000\\
-0.370000	0.630000\\
-0.360000	0.640000\\
-0.350000	0.650000\\
-0.340000	0.660000\\
-0.330000	0.670000\\
-0.320000	0.680000\\
-0.310000	0.690000\\
-0.300000	0.700000\\
-0.290000	0.710000\\
-0.280000	0.720000\\
-0.270000	0.730000\\
-0.260000	0.740000\\
-0.250000	0.750000\\
-0.240000	0.760000\\
-0.230000	0.770000\\
-0.220000	0.780000\\
-0.210000	0.790000\\
-0.200000	0.800000\\
-0.190000	0.810000\\
-0.180000	0.820000\\
-0.170000	0.830000\\
-0.160000	0.840000\\
-0.150000	0.850000\\
-0.140000	0.860000\\
-0.130000	0.870000\\
-0.120000	0.880000\\
-0.110000	0.890000\\
-0.100000	0.900000\\
-0.090000	0.910000\\
-0.080000	0.920000\\
-0.070000	0.930000\\
-0.060000	0.940000\\
-0.050000	0.950000\\
-0.040000	0.960000\\
-0.030000	0.970000\\
-0.020000	0.980000\\
-0.010000	0.990000\\
0.000000	1.000000\\
0.010000	0.990000\\
0.020000	0.980000\\
0.030000	0.970000\\
0.040000	0.960000\\
0.050000	0.950000\\
0.060000	0.940000\\
0.070000	0.930000\\
0.080000	0.920000\\
0.090000	0.910000\\
0.100000	0.900000\\
0.110000	0.890000\\
0.120000	0.880000\\
0.130000	0.870000\\
0.140000	0.860000\\
0.150000	0.850000\\
0.160000	0.840000\\
0.170000	0.830000\\
0.180000	0.820000\\
0.190000	0.810000\\
0.200000	0.800000\\
0.210000	0.790000\\
0.220000	0.780000\\
0.230000	0.770000\\
0.240000	0.760000\\
0.250000	0.750000\\
0.260000	0.740000\\
0.270000	0.730000\\
0.280000	0.720000\\
0.290000	0.710000\\
0.300000	0.700000\\
0.310000	0.690000\\
0.320000	0.680000\\
0.330000	0.670000\\
0.340000	0.660000\\
0.350000	0.650000\\
0.360000	0.640000\\
0.370000	0.630000\\
0.380000	0.620000\\
0.390000	0.610000\\
0.400000	0.600000\\
0.410000	0.590000\\
0.420000	0.580000\\
0.430000	0.570000\\
0.440000	0.560000\\
0.450000	0.550000\\
0.460000	0.540000\\
0.470000	0.530000\\
0.480000	0.520000\\
0.490000	0.510000\\
0.500000	0.500000\\
};
\addlegendentry{$\lambda_1(\VM{\omega})$}

\addplot [color=mycolor2]
  table[row sep=crcr]{%
-0.500000	1.500000\\
-0.490000	1.490000\\
-0.480000	1.480000\\
-0.470000	1.470000\\
-0.460000	1.460000\\
-0.450000	1.450000\\
-0.440000	1.440000\\
-0.430000	1.430000\\
-0.420000	1.420000\\
-0.410000	1.410000\\
-0.400000	1.400000\\
-0.390000	1.390000\\
-0.380000	1.380000\\
-0.370000	1.370000\\
-0.360000	1.360000\\
-0.350000	1.350000\\
-0.340000	1.340000\\
-0.330000	1.330000\\
-0.320000	1.320000\\
-0.310000	1.310000\\
-0.300000	1.300000\\
-0.290000	1.290000\\
-0.280000	1.280000\\
-0.270000	1.270000\\
-0.260000	1.260000\\
-0.250000	1.250000\\
-0.240000	1.240000\\
-0.230000	1.230000\\
-0.220000	1.220000\\
-0.210000	1.210000\\
-0.200000	1.200000\\
-0.190000	1.190000\\
-0.180000	1.180000\\
-0.170000	1.170000\\
-0.160000	1.160000\\
-0.150000	1.150000\\
-0.140000	1.140000\\
-0.130000	1.130000\\
-0.120000	1.120000\\
-0.110000	1.110000\\
-0.100000	1.100000\\
-0.090000	1.090000\\
-0.080000	1.080000\\
-0.070000	1.070000\\
-0.060000	1.060000\\
-0.050000	1.050000\\
-0.040000	1.040000\\
-0.030000	1.030000\\
-0.020000	1.020000\\
-0.010000	1.010000\\
0.000000	1.000000\\
0.010000	1.010000\\
0.020000	1.020000\\
0.030000	1.030000\\
0.040000	1.040000\\
0.050000	1.050000\\
0.060000	1.060000\\
0.070000	1.070000\\
0.080000	1.080000\\
0.090000	1.090000\\
0.100000	1.100000\\
0.110000	1.110000\\
0.120000	1.120000\\
0.130000	1.130000\\
0.140000	1.140000\\
0.150000	1.150000\\
0.160000	1.160000\\
0.170000	1.170000\\
0.180000	1.180000\\
0.190000	1.190000\\
0.200000	1.200000\\
0.210000	1.210000\\
0.220000	1.220000\\
0.230000	1.230000\\
0.240000	1.240000\\
0.250000	1.250000\\
0.260000	1.260000\\
0.270000	1.270000\\
0.280000	1.280000\\
0.290000	1.290000\\
0.300000	1.300000\\
0.310000	1.310000\\
0.320000	1.320000\\
0.330000	1.330000\\
0.340000	1.340000\\
0.350000	1.350000\\
0.360000	1.360000\\
0.370000	1.370000\\
0.380000	1.380000\\
0.390000	1.390000\\
0.400000	1.400000\\
0.410000	1.410000\\
0.420000	1.420000\\
0.430000	1.430000\\
0.440000	1.440000\\
0.450000	1.450000\\
0.460000	1.460000\\
0.470000	1.470000\\
0.480000	1.480000\\
0.490000	1.490000\\
0.500000	1.500000\\
};
\addlegendentry{$\lambda_2(\VM{\omega})$}

\end{axis}
\end{tikzpicture}%
		\caption{Illustration corresponding to Example \ref{ex:vb3}. If both $\VM{w}_1$ and $\VM{w}_2$ do not depend much on $\omega_1$, then $\subs{V}_1 :=\text{span}  \{ \VM{x}_1(\omega^i), i = 1,2 \} \approx \text{span} \{ \VM{x}_i(\omega^1), i = 1,2 \} =: \subs{V}_2$ but $\subs{V}_2$ is less costly to calculate than $\subs{V}_1$.} \label{fig:ex4}
	\end{center}
\end{figure}
\end{eexample}



\subsection{Upper bound of the error} \label{upperb_error}
Before we can state the algorithm, we need a measure to qualify if the subspace is rich enough or not. 
The subspace is rich enough if for every point in the parameterspace the difference between the minimal eigenvalue of the reduced problem and the minimal eigenvalue of the original large problem is smaller than a given tolerance. This means we need a way to measure the error for any given $\VM{\omega} \in \Omega$. An upper bound can be obtained only using the first eigenpair via an extension of the Bauer-Fike theorem (see \cite{Bauer1960} and more recently \cite{Saad2011}) to generalized eigenvalue problems. For this, we define the residual of an approximate eigenpair $\left( \hat{\lambda}, \VM{\hat{x}} \right)$  as $$\VM{r} = \VM{A} \VM{\hat{x}} - \hat{\lambda} \VM{B}\VM{\hat{x}}.$$
The proof of this extension is omitted as it is completely analogous to the original proof for the standard case.


%
 
\begin{theorem} \label{ste:Bauer_Fike} (Bauer-Fike for generalized eigenvalue problems) Let $\left( \tilde{\lambda}, \tilde{\VM{x}} \right)$ be an approximate eigenpair of $\left( \VM{A}, \VM{B}\right)$ where $\tilde{\VM{x}}$ is of $\VM{B}$-norm unity and $\tilde{\lambda} = \left( \VM{A} \tilde{\VM{x}}, \tilde{\VM{x}} \right)$. Then, there exists an eigenvalue $\lambda$ of $\left( \VM{A}, \VM{B} \right)$ such that
\begin{equation} \label{eqn:Ext_Bauer_Fike} \mid { \tilde{\lambda} - \lambda} \mid \leq  \dfrac{ \left \| \VM{r} \right \|_2 }{ \sqrt{\lambda_1(\VM{B})}}.
\end{equation}
\end{theorem}
The residual norm is very cheap to calculate for a set of parameter samples if the dependency of $\VM{A}$ and $\VM{B}$ is affine. More precisely, if $\VM{A}(\VM{\omega})$ and $\VM{B}(\VM{\omega})$ can be written as $\VM{A}(\VM{\omega}) = \sum_{i=1}^{m_0} \theta_{A,i}(\VM{\omega}) \VM{A}_i$ resp. $\VM{B}(\VM{\omega}) = \sum_{i=1}^{m_1} \theta_{B,i}(\VM{\omega}) \VM{B}_i$ with $\theta_{A,i}, \theta_{B,j}, i= 1, \hdots m_0, j = 1, \hdots, m_1$ scalar functions then the residual is in this case
\begin{align*}
\VM{r}(\VM{\omega}) & = \VM{A}(\VM{\omega}) \tilde{\VM{x}}(\VM{\omega}) - \tilde{\lambda}(\VM{\omega}) \VM{B}(\VM{\omega})  \tilde{\VM{x}}(\VM{\omega}) \\
& = \sum_{i=1}^{m_0} \theta_{A,i}(\VM{\omega}) \left( \VM{A}_i \VM{V}\right) \VM{x}^\subs{V}(\VM{\omega}) - \tilde{\lambda}(\VM{\omega}) \sum_{i=1}^{m_1} \theta_{B,i}(\VM{\omega}) \left( \VM{B}_i \VM{V}\right) \VM{x}^\subs{V}(\VM{\omega}).
\end{align*}
This means that the precomputation of $\VM{A}_i, i = 1, \hdots, m_0$ and $\VM{B}_i, i = 1, \hdots, m_1$ reduces the cost of calculating the residual norm.
The drawback of \eqref{eqn:Ext_Bauer_Fike} is that this is usually a very crude bound. If we have also information about the second minimal eigenvalue, we can use the Kato-Temple theorem \cite{Kato1949}, \cite{Temple1952} for a sharper bound. We rephrase the theorem for the generalised eigenvalue problem.

\begin{theorem} \label{ste:KatoTemple}
(Kato-Temple for generalized eigenvalue problems)  Let $\left( \tilde{\lambda}, \tilde{\VM{x}} \right)$ be an approximate eigenpair of $\left( \VM{A}, \VM{B}\right)$  where $\tilde{\VM{x}}$ is of $\VM{B}$-norm unity and $\tilde{\lambda} = \left( \VM{A} \tilde{\VM{x}}, \tilde{\VM{x}} \right)$. Assume that we know an interval $]\alpha, \beta[$ that contains $\tilde{\lambda}$ and one eigenvalue $\lambda$ of $\left( \VM{A}, \VM{B} \right)$. Then it holds
$$ - \dfrac{ \norm{\VM{r}}_2^2}{\lambda_1(\VM{B}) \left( \tilde{\lambda} - \alpha \right)} \leq \tilde{\lambda}-\lambda \leq \dfrac{ \norm{\VM{r}}_2^2}{\lambda_1(\VM{B})\left( \beta - \tilde{\lambda} \right)}.$$
\end{theorem}

If we now define $\delta$ as the distance to the second minimal eigenvalue, then we get 
\begin{equation} |\tilde{\lambda} - \lambda_1 | \leq \dfrac{ \norm{\VM{r}}_2^2}{\lambda_1(\VM{B})\delta},\label{eqn:KatoTemple} \end{equation}
which is in general a sharper upper bound then the one from Theorem \ref{ste:Bauer_Fike}.

In the Bauer-Fike as well as in the Kato-Temple theorem the bound depends on the first eigenvalue of $\VM{B}$. As it is expensive to calculate this eigenvalue for every $\VM{\omega}$, we only calculate it for one parametervalue and use it for the whole domain. As already mentioned, we approximate the second, respectively third eigenvalue by the second, respectively third Ritz value from the Krylov subspace. Whether we use the Bauer-Fike theorem or the Temple-Kato theorem depends on how accurately we can determine the gap $\delta$ between the first two eigenvalues. If we add the first m eigenvectors $(m > 1)$ corresponding to every sample point to the subspace, we interpolate the first m eigenvalues in the sample points, and obtain in this way a global approximation of the first m eigenvalues over the parameter space.  Hence we expect that the gap $\delta$ can be estimated well from the projected eigenvalue problem,  and we make in this case use of the Kato-Temple theorem, where we replace $\delta$ by its estimate. If only one eigenvector is added in a sample point, an accurate estimate of the gap from the projected eigenvalue problem cannot be guaranteed. In the latter case the Bauer-Fike theorem is invoked instead.

\subsection{Algorithm} \label{sect:algo}
We now have all ingredients for an algorithm to compute the reduced problem. A remaining crucial element is the choice of sample points, which we now discuss. The proposed method is inspired by the reduced basis method for solving parametric
 PDEs, see \cite{Quarteroni2016}. 

We assume that the parameter space $\Omega$ is a Cartesian product of intervals in $\R$, so $$\Omega = \Omega_1 \times \Omega_2 \times \hdots \times \Omega_d = [a^1, b^1] \times [a^2, b^2] \times \hdots \times [a^d, b^d].$$ We first make an initial subspace $\subs{V}$ constructed with the first eigenvector(s) and the partial derivatives for all sample points in an initial set $\Omega_\text{init}$ (see line 1-4 in Algorithm \ref{algo:Redbasis}) and we check then if this constructed subspace is large enough by testing the upper bound on some training set $\Omega_\text{train}$. Whether we use the Kato-Temple theorem \ref{ste:KatoTemple} or the Bauer-Fike theorem \ref{ste:Bauer_Fike} depends on the number of eigenvectors we add per interpolation point, see previous section.
 If the upper bound at all points in the training set $\Omega_\text{train}$ is not below a certain tolerance, we add the first eigenvectors and partial derivatives of the point in the training set where the upper bound is maximal (see line 13-22 in Algorithm \ref{algo:Redbasis}). 
Once the upper bound is below the tolerance, these sample points are removed from the training set  (see line 7-10 in Algorithm \ref{algo:Redbasis}). We are allowed to do this because this implies that the error on the eigenvalue is below the tolerance and furthermore the error cannot increase by the second assertion in Theorem \ref{prop:interpol} . 
The initial set of sample points is decomposed as 
\begin{equation}  \label{eqn:algo_init2}
\Omega_\text{init} := \Omega^\text{init}_1 \times \Omega^\text{init}_2 \times \hdots \times \Omega^\text{init}_d
\end{equation} 
where \begin{equation} \label{eqn:algo_init1}
\Omega_i^\text{init} = \{ a_i + (j-1) h_i | j = 1, 2, \hdots, n_i\}, h_i = \dfrac{b^i - a^i}{n_i-1}, i = 1, 2, \hdots, d.
\end{equation}
This is an initial set of $n_\text{init} = n_1 n_2 \hdots n_d$ points. The training set $\Omega_\text{train}$ is composed in a similar way but using a much finer grid. The algorithm is stated in pseudocode in Algorithm \ref{algo:Redbasis}. As we do experiments with adding partial derivatives to the subspace and with not adding them to the subspace, we state 'optionally' after the lines where we calculate the partial derivatives. 

\begin{algorithm}[h] 
\caption{Aim: Calculating a subspace$ \subs{V}$ such that $ | \lambda^\subs{V}_1(\VM{\omega}) - \lambda_1(\VM{\omega})| < \text{tol}, \forall \VM{\omega} \in \Omega$} \label{algo:Redbasis}
\hspace*{\algorithmicindent} \textbf{Input:}  \begin{enumerate}
\item Matrix $\VM{A}(\VM{\omega})$ and $\VM{B}(\VM{\omega})$
\item All partial derivatives $\pderiv{\VM{A}(\VM{\omega})}{\omega_i}$ and $\pderiv{\VM{B}(\VM{\omega})}{\omega_i}, i = 1,2,\hdots,d$
\item Choose an initial set $\Omega_\text{init}$ and a training set $\Omega_\text{train}$ of sample points.
\item $n_\text{max}$ maximal number of iterations
\end{enumerate}
\hspace*{\algorithmicindent} \textbf{Output:} $\lambda^\subs{V}_1(\VM{\omega})$ such that $\max_{\VM{\omega} \in \Omega_\text{train}} | \lambda^\subs{V}_1(\VM{\omega}) - \lambda_1(\VM{\omega})| < \text{tol} $
\begin{algorithmic}[1]
\State Calculate $\left( \lambda_1(\VM{\omega}), \VM{x}_1(\VM{\omega}) \right), \VM{\omega} \in \Omega_\text{init}$ 
\State Estimate $\left( \tilde{\lambda}_j(\VM{\omega}), \tilde{\VM{x}}_j(\VM{\omega}) \right), j = 2, 3, \hdots, m, \VM{\omega} \in \Omega_\text{init}$ 
\State Calculate partial derivatives $\pderiv{\VM{x}_1(\VM{\omega})}{\omega_j}, j = 1, 2, \hdots, d, \VM{\omega} \in \Omega_\text{init}$ from system \eqref{eqn:system_dx} if $|\lambda_1 - \tilde{\lambda}_2|> 10^{-8}$  (If $\lambda_1$ is simple) (Optionally)
\State $ \subs{V}:= \text{span}\{ \VM{x}_1(\VM{\omega}), \tilde{\VM{x}}_2(\VM{\omega}) \hdots, \tilde{\VM{x}}_m(\VM{\omega}),  \pderiv{\VM{x}_1(\VM{\omega})}{\omega_1}, \pderiv{\VM{x}_1(\VM{\omega})}{\omega_2}, \hdots, \pderiv{\VM{x}_1(\VM{\omega})}{\omega_d} | \VM{\omega} \in \Omega_\text{init} \} $
\State $u_\text{old}(\VM{\omega}) = 1, \VM{\omega} \in \Omega_\text{train}$
\For{$i=1, 2, \hdots, n_\text{max}$}
\State \underline{Update the training set}
\State Calculate new upper bound $u_\text{new}(\VM{\omega})$ for all $\VM{\omega} \in \Omega_\text{train} $ using Theorem \ref{ste:Bauer_Fike} if $m = 1$ or Theorem \ref{ste:KatoTemple} if $m > 1$
\State Upper bound $u(\VM{\omega}) = \min( u_\text{old}(\VM{\omega}), u_\text{new}(\VM{\omega}))$
\State $\Omega_\text{valid} = \{ \VM{\omega} | \VM{\omega} \in \Omega_\text{train} \text{ and } u(\VM{\omega}) < \text{tol} \}$
\State $\Omega_\text{train} = \Omega_\text{train} \setminus \Omega_\text{valid}$
\If{$\Omega_\text{train}$ empty}
\State break
\Else
\State \underline{Adding vectors to the subspace}
\State $\VM{\omega}^i = \text{arg} \max_{\VM{\omega} \in \Omega_\text{train}} u(\VM{\omega}) $
\State Calculate $\left( \lambda_1(\VM{\omega}^i), \VM{x}_1(\VM{\omega}^i) \right)$
\State Estimate $\left( \tilde{\lambda}_j(\VM{\omega}^i), \tilde{\VM{x}}_j(\VM{\omega}^i \right), j = 2,3, \hdots, m$
\If{$|\lambda_1 - \tilde{\lambda}_2|> 10^{-8}$}  (Check if $\lambda_1$ is simple)
\State Calculate partial derivatives $\pderiv{\VM{x}_1(\VM{\omega}^i)}{\omega_j}, j = 1,2, \hdots, d$ from system \eqref{eqn:system_dx} (optionally) \\
\State $\subs{V} = \subs{V} \bigcup \{ \VM{x}_1(\VM{\omega}^i), \tilde{\VM{x}}_2(\VM{\omega}^i) \hdots, \tilde{\VM{x}}_m(\VM{\omega}^i), \} \left( \bigcup \{  \pderiv{\VM{x}_1(\VM{\omega}^i)}{\omega_1}, \pderiv{\VM{x}_1(\VM{\omega}^i)}{\omega_2}, \hdots, \pderiv{\VM{x}_1(\VM{\omega}^i)}{\omega_d} \} \text{(optionally)} \right)$
\Else
\State $\subs{V} = \subs{V} \bigcup \{ \VM{x}_1(\VM{\omega}^i), \tilde{\VM{x}}_2(\VM{\omega}^i) \hdots, \tilde{\VM{x}}_m(\VM{\omega}^i) \}$
\EndIf
\State Make basis for $\subs{V}$
\State $\Omega_\text{train} = \Omega_\text{train} \setminus \VM{\omega}^i$
\EndIf 
\EndFor
 \end{algorithmic}
\end{algorithm}

\subsection{Saturation assumption}
Although the training set becomes smaller in every iteration, calculating the upper bound for all sample points in the training set can become expensive if many iterations are needed. Instead of updating all upper bounds, we select only those that potentially generate the highest upper bound for the error.
We therefore make use of the so-called saturation assumption, see \cite{Hesthaven2014}. Let $u^k$ be the upper bound at iteration $k$. The saturation assumption says that there exists a $C > 0$ such that 
\begin{equation} \label{eqn:satur}
u^l(\VM{\omega}) < C u^k(\VM{\omega}), \forall l > k, \forall \VM{\omega} \in \Omega.
\end{equation}
In our method the saturation assumption is fulfilled for $C = 1$ as we only update the upper bound $u^k$ if the upper bound has decreased. We use this in the following way.  
We first sort the parametervalues in $\Omega_\text{train}$  by the upper bound of the previous iteration in descending order. We recompute the upper bound and update after each computation the maximal upper bound $u_\text{max}$. If for a certain $\VM{\omega}$, it holds that $u^k(\VM{\omega}) < u_\text{max}$ we are allowed to skip all the next sample points in the sequence by \eqref{eqn:satur}. The next $\VM{\omega}$ for which we add the eigenvectors and all partial derivatives, is the $\VM{\omega}$ which has the upper bound $u_\text{max}$. 

\section{Implementation details and numerical results}
All algorithms written in the previous sections are implemented in Matlab version R2017a. All experiments are performed on an Intel i5-6300U with 2.5 GHz and 8 GB RAM. In all examples, the linear systems in building the Krylov space and the systems of the form \eqref{eqn:system_dx} to compute the derivative of an eigenvector could still be solved efficiently using a direct method. Practically, because the used matrices are sparse, it is efficient to use the backslash-operator in Matlab for solving the systems of the form \eqref{eqn:system_dx}. We use the \texttt{eigs}-command in Matlab to calculate the minimal eigenvalue which uses the shift-and-invert Arnoldi's method with shift 0. This method makes a Krylov subspace and we adapt the code such that it also returns this space. We use the other vectors in the Krylov subspace to approximate the second eigenvector. In all examples a tolerance for the error of $10^{-5}$ is used. 
We applied our algorithm to examples where the minimal eigenvalue corresponds to calculating the coercivity constant when solving PDEs \cite{Sirkovic2016} and to one example from structural mechanics. We compare the results from four cases: adding only the first eigenvector, adding the first two eigenvectors, adding the first eigenvector and all partial derivatives and finally the first two eigenvectors and all partial derivatives. In the next tables, we put the dimension of the found subspace, the number of sample points, the total computing time and the average time we needed to compute the partial derivatives and eigenvectors. \\

For Example \ref{vbd:2}, the initial set is built by discretising both dimensions of the parameterspace into $4$ to get an initial subset of $16$ points. For the training set, we discretise the first respectively the second dimension of the parameterspace into $25$ respectively $40$ points, to have a total of $1000$ points in our training set. The results can be found in Table \ref{tbl:vbd2}. We see that adding other information than only the first eigenvector is beneficial. The number of points where we need to calculate the second eigenvector and the partial derivative is much lower which results in a smaller computational time. The main benefit from adding the second eigenvector is that we can use the more strict upper bound of Kato-Temple. We also noticed that the time needed to compute the derivatives is much lower than calculating the eigenvalues. We see we get the best results when we add both the eigenvectors and the partial derivatives into the subspace.

\begin{table}[h] 
	\centering
\caption{The results for Example \ref{vbd:2} \label{tbl:vbd2}}
\begin{tabular*}{500pt}{@{\extracolsep\fill}lcccc@{\extracolsep\fill}}
\toprule
&1 eigv &2 eigv &1 eigv + partial deriv. &2 eigv + partial deriv. \\
\midrule
dimension $\subs{V}$&96&108&138&132\\
nbr points&96&54&46&33\\
total time&59.040s&38.923s&49.997s&29.293s\\
time derivative (per vector)&/&/&0.002s&0.002s\\
time eigenv (per vector)&0.080s&0.042s&0.079s&0.044s\\
\bottomrule
\end{tabular*}
 
\end{table}

The next example is also taken from \cite[Example 4.3]{Sirkovic2016}. The parameterspace is $\Omega = [-0.1, 0.1] \times [0.2, 0.3]$ and the matrix $\VM{A}$ is of the form $\VM{A}(\VM{\omega}) = \sum_{i=1}^{16} \theta_i(\VM{\omega}) \VM{A}_i$ with $\theta_i(\VM{\omega})$ analytic functions and matrix $\VM{B}$ is constant. The dimension of the problem is $2183$. As training set we choose a grid where we discretise the interval for $\omega_1$ into $40$ and $\omega_2$ into $25$ points. As initial set we discretise both intervals into 3. 
The results are stated in Table \ref{tbl:coerc_1}. We see that in this case the needed subspace is low-dimensional which results in a small computational time. Similar conclusions can be drawn for this example.
\begin{table}[h] 
	\centering
\caption{Results \label{tbl:coerc_1}}
\begin{tabular*}{500pt}{@{\extracolsep\fill}lcccc@{\extracolsep\fill}}
\toprule
&1 eigv &2 eigv &1 eigv + partial deriv. &2 eigv + partial deriv. \\
\midrule
dimension $\subs{V}$&37&30&48&40\\
nbr points&37&15&16&10\\
total time&32.038s&14.617s&24.372s&8.499s\\
time derivative (per vector)&/ & /&0.005s&0.004s\\
time eigenv (per vector)&0.034s&0.015s&0.026s&0.012s\\
\bottomrule
\end{tabular*}

\end{table} 

\begin{figure}[h]
\begin{center}
\setlength{\figW}{5cm} 
  \input{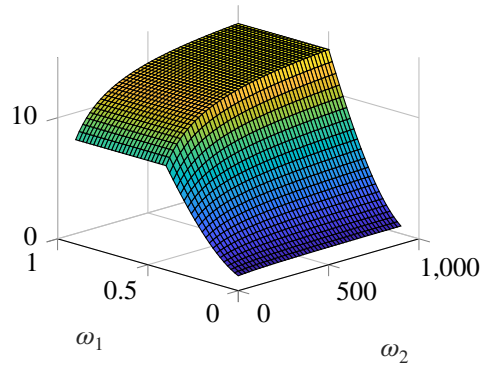}
\caption{The minimal eigenvalue $\lambda_1(\VM{\omega})$ over the parameter for the example from structural mechanics. We see the minimal eigenvalue does not smoothly depend on the variable. } \label{vb:bouwkunde}
\end{center} 
\end{figure}

In the previous examples matrix $\VM{B}$ was constant, since the article from which we took the examples considered only eigenvalue problems with constant $\VM{B}$. The next example originates from an application in structural mechanics and here matrix $\VM{B}$ depends on $\VM{\omega}$. \\
The matrices in this example are the stiffness and mass matrices arising from a finite element approximation of a vibrating beam that can deflect in both directions perpendicular to its own axis.  The first parameter concerns the size of the beam along one of the axis (in interval $[0.1,1]$) and the second parameter deals with the stiffness and density of the material (in interval $[100, 1000]$). For the initial set, we discretise both dimensions into three. The matrices have size $1404$ and they both depend on at least one of the parameters. In this problem ,the stiffness and mass matrix both dependent on the two parameters. One is usually interested in the minimal eigenfrequencies which boils down to calculating the minimal eigenvalue. In Figure \ref{vb:bouwkunde} we depict the minimal eigenvalue over the parameterspace and Table \ref{tbl:bouwkunde} summarizes the results. We observe that only a small subspace was needed and that we benefit from including partial derivatives in the subspace since the computation of the eigenvectors takes much more time than solving the system to compute the partial derivatives.

\begin{table}[h] 
	\centering
\caption{Results for example from structural mechanics \label{tbl:bouwkunde}}
\begin{tabular*}{500pt}{@{\extracolsep\fill}lcccc@{\extracolsep\fill}}
\toprule
&1 eigv &2 eigv &1 eigv + partial deriv. &2 eigv + partial deriv. \\
\midrule
dimension $\subs{V}$&17&  22&32& 31\\
nbr points&17&11&12&9\\
total time&193.727s&105.806s&123.125s&83.531s\\
time derivative (per vector)&/&/&0.005s&0.005s\\
time eigenv (per vector)&1.743s&0.871s& 1.787s&0.890s\\
\bottomrule
\end{tabular*}

\end{table} 

\section{Conclusion}
The main contribution of this paper is that we showed that for calculating extreme eigenvalues of eigenpairs $(\VM{A}, \VM{B})$ where both matrices are symmetric and $\VM{B}$ is positive definite, it is beneficial to add, besides the first eigenvector, also its partial derivatives, or alternatively, the Ritz vector associated with the second and third Ritz values. We have proved that in the case of the addition of partial derivatives, the eigenvalues satisfy a Hermite interpolation property of order 2. Numerical examples confirm this statement and showed that this property can also be used in practice.

The question arises what we do if the system of linear equations for calculating the partial derivatives cannot be efficiently solved by direct solvers. A possibility is to first estimate the partial derivative from the Krylov subspace $\subs{V}_K$ obtained when calculating the first eigenvector as this space contains approximations of the vectors from which the Ritz-values are the closest to the minimal eigenvalue. Let $\VM{V}_K$ be a basis of this subspace then we first project the system in \eqref{eqn:system_dx} on $\subs{V}_K$ to obtain

\begin{equation} \label{eqn:init_sol_GMRES}
\begin{bmatrix}
\VM{V}_K^T \left( \lambda_1  \VM{B}-\VM{A} \right) \VM{V}_K & \VM{V}_K^T \VM{B} \VM{x}_1 \\ 
\VM{x}_1^T \VM{B} \VM{V}_K & 0
\end{bmatrix} \begin{bmatrix} \pderiv{\VM{x}_1}{\omega_i}^{\subs{V}_K} \\ \pderiv{\lambda^{\subs{V}_K}_1}{\omega_i} \end{bmatrix} = \begin{bmatrix}
\VM{V}_K^T \left(\pderiv{\VM{A}}{\omega_i}   - \lambda_1 \pderiv{\VM{B}}{\omega_i}  \right ) \VM{x}_1 \\
- \dfrac{ \VM{x}_1^T \pderiv{\VM{B}}{\omega_i} \VM{x}_1}{2} 
\end{bmatrix}
\end{equation}
which has the same dimension as $\subs{V}_K$.
We can use this approximation as a starting vector in an iterative method like GMRES to compute the derivative.

\section*{Acknowledgements}
The authors thank the referees for their helpful remarks.
This work was supported by the project KU Leuven Research Council grant OT/14/074 and C14/17/072 and by the project G0A5317N of the Research Foundation-Flanders (FWO - Vlaanderen).
\bibliography{references}%

%

\end{document}